\documentclass[12pt,a4paper]{article}

\usepackage{amsthm}
\usepackage{epsfig}
\usepackage{alltt}
\usepackage{makeidx}
\usepackage{newlfont}
\usepackage{amsmath}
\usepackage{amssymb}
\usepackage{amsfonts}
\usepackage{amscd}
\usepackage[ngerman]{babel}
\usepackage{hyperref}
\usepackage{verbatim}
\usepackage{changebar}
\input{xy}
\xyoption{all}

\makeindex

\newtheorem{Satz}{Satz}[subsection]

\newtheorem{Proposition}[Satz]{Proposition}
\newtheorem{Korollar}[Satz]{Korollar}
\newtheorem{Lemma}[Satz]{Lemma}
\newtheorem{Theorem}[Satz]{Theorem}

\theoremstyle{definition}
\newtheorem{Definition}[Satz]{Definition}

\theoremstyle{remark}
\newtheorem{Bemerkung}[Satz]{Bemerkung}

\newtheorem{Notation}[Satz]{Notation}

\newcommand{\op}{\operatorname}

\newcommand{\ra}{\rightarrow}

\newcommand{\lra}{\longrightarrow}
\newcommand{\sra}{\twoheadrightarrow}
\newcommand{\hra}{\hookrightarrow}

\newcommand{\sira}{\stackrel{\sim}{\rightarrow}}

\newcommand{\da}{\downarrow}

\newcommand{\IFF}{\Leftrightarrow}

\newcommand{\Bl}[1]{{\mathbb{#1}}}

\newcommand{\DC}{\Bl{C}}

\newcommand{\defind}[1]{{\bf #1}\index{#1}}

\newcommand{\As}{\;{\mbox {$\subset \hspace{-2,3ex}$  \raisebox{-0,5ex}
{${\large A}$} }}}

\xyoption{all}
\begin{document}
\title{Andersen-Filtrierung und harter Lefschetz}
\author{Wolfgang Soergel}

\maketitle 
\begin{center}
  \emph{F"ur Joseph Bernstein}
  \end{center}

\begin{abstract}
We consider the principal block of category $\cal{O}$ and its $\Bbb{Z}$-graded
version introduced in \cite{BGSo}.
On the space of homomorphisms from a Verma module to an indecomposable tilting
module we may define natural filtrations following Andersen \cite{AFil}.
The arguments given in this article prove that these 
filtrations are compatible with the $\Bbb{Z}$-graded
structure, although explicitely we only show that the
dimensions of the sucessive subquotients of the
filtration are compatible with this intuition. 
This statement is very similar to the 
semisimplicity of the subquotients of
the Jantzen filtration proved in \cite{BB-J}, 
but the method of proof is quite different.
I would like to know how to directly relate both results,
as this would give an alternative proof of said semisimplicity.
\end{abstract}
\section{Deformationen der Kategorie $\cal{O}$}
\begin{Bemerkung}
  In diesem und dem n"achsten Abschnitt wiederholen wir Resultate von \cite{GJ}
  in einer f"ur unsere Ziele angepa"sten Sprache,
die auch \cite{FieC} sehr nahe steht.  Seien $\frak{g} \supset
  \frak{b}\supset \frak{h}$ eine halbeinfache komplexe Lie-Algebra, eine
  Borel'sche und eine Cartan'sche.  Sei $S = S \frak{h} 
= \cal{O}( \frak{h}^{\ast})$ die symmetrische Algebra von $\frak{h}$.  
Wir betrachten die Kategorie $\op{Kring}^S$ aller kommutativen
unit"aren Ringe $T$ mitsamt einem ausgezeichneten Homomorphismus
$\varphi : S \ra T.$ Gegeben $T\in \op{Kring}^S$
  betrachten wir die Kategorie $\frak{g}\op{-Mod_{\DC}-} T$ aller
  $\frak{g}\op{-}T$-Bimoduln, auf denen die Rechts- und die Linksoperation von
  $\Bbb{C}$ zusammenfallen.  Nicht weiter spezifizierte Tensorprodukte 
sind stets "uber
  $\DC$ zu verstehen.
\end{Bemerkung}
\begin{Definition}
Gegeben $T=(T,\varphi)\in\op{Kring}^S$  
erkl"aren f"ur alle $\lambda \in \frak{h}^{\ast}$ in 
jedem Bimodul $M \in
  \frak{g}\op{-Mod_{\DC}-}T$ den {\bf deformierten Gewichtsraum}
  $M^{\lambda}$  durch die Vorschrift
  $$M^{\lambda} =M^{\lambda}_T=\{ m \in M 
\mid (H - \lambda (H)) m = m \varphi (H)
  \quad\forall H \in \frak{h}\}$$
\end{Definition}
\begin{Bemerkung} 
Gegeben $M \in \frak{g}\op{-Mod_{\DC}-}T $ ist die 
kanonische Abbildung von der direkten Summe seiner
deformierten Gewichtsr"aume  nach $M$ stets 
eine Injektion $\bigoplus_\lambda M^\lambda\hra M.$
Im Fall $T = \cal{O}( \frak{h}^{\ast})$ ist das offensichtlich, da die 
Gewichtsr"aume $M^{\lambda}$
als $T \otimes T$-Moduln gerade Tr"ager im Graphen von 
$(\lambda +): 
\frak{h}^{\ast} \ra \frak{h}^{\ast}$ haben und da diese Graphen paarweise 
disjunkt sind.
Im allgemeinen haben unsere Gewichtsr"aume Tr"ager im Urbild unserer 
Graphen unter der von
$\op{id} \otimes \varphi$ induzierten Abbildung 
$\op{Spec} (\Bbb{C} [\frak{h}^{\ast}] \otimes
T) \ra \op{Spec} (\cal{O}( \frak{h}^{\ast})
\otimes\cal{O}( \frak{h}^{\ast}))$
und sind damit ebenfalls disjunkt.
\end{Bemerkung}
\begin{Definition}
 F"ur jedes $T\in\op{Kring}^S$  definieren wir 
in unserer Kategorie von Bimoduln   eine volle Unterkategorie,
  die {\bf deformierte Kategorie}
  $$\cal{O}(T) \subset \frak{g}\op{-Mod_{\DC}-}T$$
  als die Kategorie aller
  Bimoduln, die lokal endlich sind unter $\frak{n}  = [\frak{b},\frak{b}]$
  und die die Summe ihrer deformierten Gewichtsr"aume sind.
\end{Definition}
\begin{Bemerkung}  
Besonders prominente Objekte dieser Kategorie sind die {\bf deformierten
  Vermamoduln}
  $$\Delta_{T} (\lambda) =
  \op{prod}_{\frak{b}}^\frak{g}(\Bbb{C}_{\lambda}\otimes T) =U (\frak{g})
  \otimes_{U(\frak{b})} (\Bbb{C}_{\lambda}\otimes T)$$
  f"ur $\lambda \in
  \frak{h}^{\ast}$, wobei 
zu verstehen ist, da"s die Rechtsoperation von $T$ nur
  den letzten Tensorfaktor bewegt, die Operation von $U(\frak{b})$ auf
  $\Bbb{C}_{\lambda} \otimes T$ aber vermittels der kanonischen Surjektion
  $\frak{b} \twoheadrightarrow \frak{h}$ von der Tensoroperation von $\frak{h}$
  herkommt, wobei $H \in \frak{h}$ auf $\Bbb{C}_{\lambda}$ 
operieren m"oge durch
  den Skalar $\lambda (H)$ und auf $T$ durch Multiplikation mit $\varphi (H)$.
\end{Bemerkung}
\begin{Bemerkung}
  Die Kategorie $\cal{O} (T)$ ist stabil unter dem Tensorieren von links mit
  endlichdimensionalen Darstellungen von $\frak{g},$ wobei wir als
  Linksoperation von $\frak{g}$ auf solch
einem Tensorprodukt die Tensoroperation nehmen und als Rechtsoperation
  von $T$  die Rechtsoperation auf dem zweiten Tensorfaktor.  Mit
  einem Bimodul enth"alt $\cal{O} (T)$ auch alle seine Subquotienten.  Im
  Spezialfall $T = \Bbb{C}$ und $\varphi$ dem Auswerten am Nullpunkt von $
  \frak{h}^{\ast}$ ist $\cal{O} (T),$ wenn man vom Fehlen gewisser
  Endlichkeitsbedingungen einmal absieht, schlicht die "ubliche Kategorie
  $\cal{O}$ von Bernstein-Gelfand-Gelfand, und
  $\Delta_\DC(\lambda)=\Delta(\lambda)$ ist der 
Vermamodul mit h"ochstem Gewicht
  $\lambda.$
\end{Bemerkung}
\begin{Definition}
  Wir betrachten weiter die bez"uglich $\frak{h}$ zu $\frak{b}$ opponierte
  Borel'sche $\bar{\frak{b}} \subset \frak{g}$ und 
f"ur $\lambda\in \frak{h}^\ast$ den Unterbimodul
  $$\nabla_{T} (\lambda) \subset
  \op{ind}_{\bar{\frak{b}}}^\frak{g}(\Bbb{C}_{\lambda}\otimes T)
  =\op{Hom}_{U(\bar{\frak{b}})} (U (\frak{g}), \Bbb{C}_{\lambda}\otimes T)$$
  der
  als die Summe aller deformierten Gewichtsr"aume im fraglichen $\op{Hom}$-Raum
  erkl"art wird.  Wir nennen ihn den {\bf deformierten Nabla-Modul} 
mit h"ochstem
  Gewicht $\lambda.$
\end{Definition}
\begin{Bemerkung}
  Unter der durch die Restriktion gegebenen Identifikation unseres
  $\op{Hom}$-Raums mit $\op{Hom}_{\Bbb{C}}(U (\frak{n}), \Bbb{C}_{\lambda}
  \otimes T)$ entspricht $\nabla_{T}(\lambda)$ genau denjenigen Homomorphismen,
  die nur auf endlich vielen $\frak{h}$-Gewichtsr"aumen aus $U(\frak{n})$
  von Null verschieden sind. Die deformierten Nabla-Moduln geh"oren auch zu
  $\cal{O} (T).$
\end{Bemerkung}

\begin{Bemerkung}
  Alle Gewichtsr"aume von $\nabla_{T} (\lambda)$ 
und $\Delta_{T} (\lambda)$ sind
  "uber $T$ frei und endlich erzeugt, und ist $T$ nicht der Nullring, so haben
  die deformierten Gewichtsr"aume zum Gewicht $(\lambda - \nu)$ 
in beiden Moduln
  den Rang $\op{dim}_{\Bbb{C}}U(\frak{n} )^{\nu}.$ Wir haben kanonische
  $T$-Modul-Morphismen $T \overset{\sim}{\ra} \Delta_{T}(\lambda)^{\lambda}
  \hookrightarrow \Delta_{T} (\lambda)$ sowie $\nabla_{T}(\lambda)
  \twoheadrightarrow \nabla_{T}(\lambda)^{\lambda} \overset{\sim}{\ra} T$ und
  f"ur jede Ringerweiterung $T \ra T^{\prime}$ kanonische Isomorphismen
  $\Delta_{T} (\lambda) \otimes_{T} T^{\prime} \overset{\sim}{\ra}
  \Delta_{T^{\prime}} (\lambda)$ sowie $\nabla_{T} (\lambda) \otimes_{T}
  T^{\prime} \overset{\sim}{\ra} \nabla_{T^{\prime}} (\lambda).$
\end{Bemerkung}

\begin{Bemerkung}
  Wir w"ahlen nun f"ur unsere Lie-Algebra einen involutiven Automorphismus 
$\tau: \frak{g} \ra \frak{g}$
  mit der Eigenschaft $\tau |_{\frak{h}} = -\op{id}$ und definieren einen
  kontravarianten Funktor
  $$d= d_{\tau} : \frak{g}\op{-Mod_{\DC}-}T \ra \frak{g}\op{-Mod_{\DC}-}T$$
  durch die Vorschrift, da"s $dM \subset \op{Hom}_{-T} (M,T)^{\tau}$ die Summe
  aller deformierten Gewichtsr"aume sein soll im fraglichen $\op{Hom}$-Raum mit
  der durch $\tau$ vertwisteten $\frak{g}$-Operation.  Ist $M\in
  \frak{g}\op{-Mod_{\DC}-}T $ die Summe seiner deformierten
  Gewichtsr"aume, so haben wir einen kanonischen 
Morphismus $M\ra ddM,$ und sind
  zus"atzlich alle deformierten Gewichtsr"aume von $M$ frei und endlich erzeugt
  "uber $T,$ so ist dieser Morphismus ein Isomorphismus.
\end{Bemerkung}

\begin{Bemerkung}
  Die Restriktion auf den h"ochsten 
deformierten Gewichtsraum definiert zusammen
  mit der universellen Eigenschaft der 
induzierten Darstellung einen kanonischen
  Homomorphismus
  $$
  \op{Hom}_{-T} (\op{prod}^{\frak{g}}_{\bar{\frak{b}}} (\Bbb{C}_{\lambda}
  \otimes T),T) \sira \op{ind}_{\bar{\frak{b}}}^{\frak{g}} \op{Hom}_{-T}
  (\Bbb{C}_{\lambda} \otimes T,T)
  $$
  von dem man durch Betrachtung der deformierten Gewichtsr"aume erkennt, da"s
  er einen Isomorphismus von Bimoduln $$d\Delta_{T} (\lambda) \cong \nabla_{T}
  (\lambda)$$
  liefert. Mit unseren Vor"uberlegungen folgt dann auch $d\nabla_{T}
  (\lambda) \cong \Delta_{T} (\lambda) .$ Nach der Tensoridentit"at hat weiter
  $E \otimes \Delta_{T} (\lambda)$ eine Filtrierung mit Subquotienten
  $\Delta_{T} (\lambda + \nu),$ wo $\nu$ "uber die Multimenge $ P (E)$ der
  Gewichte von $E$ l"auft, und da $E \otimes $ bis auf die
Wahl eines Isomorphismus $dE\cong  E$ mit 
der Dualit"at $d$ vertauscht,
  gilt Analoges f"ur $E\otimes  \nabla_{T} (\lambda).$
\end{Bemerkung}

\begin{Proposition}\label{EHo}
  \begin{enumerate}
  \item F"ur alle $\lambda$ induziert die Restriktion auf den deformierten
    Gewichtsraum zum Gewicht $\lambda$ zusammen mit den beiden kanonischen
    Identifikationen $\Delta_{T}(\lambda)^{\lambda} \overset{\sim}{\ra} T$ und
    $\nabla_{T} (\lambda)^{\lambda} \overset{\sim}{\ra}T$ einen Isomorphismus
    $$\op{Hom}_{\cal{O}(T)} (\Delta_{T} (\lambda),
      \nabla_{T}(\lambda)) \overset{\sim}{\ra} T$$
    \item F"ur $\lambda \neq \mu$ aus $\frak{h}^{\ast}$ 
          gilt
      $\op{Hom}_{\cal{O} (T)} (\Delta_{T} (\lambda),
      \nabla_{T} (\mu))=0.$
\item F"ur alle $ \lambda, \mu \in\frak{h}^{\ast}$ gilt
$\op{Ext}^{1}_{\cal{O} (T)} (\Delta_{T} (\lambda), \nabla_{T} (\mu))=0 .$ 
\end{enumerate}
\end{Proposition}
\begin{proof}[Beweis]
Bezeichne $R^+\subset \frak{h}^\ast$ die Wuzeln von $\frak{n} $ 
und $|R^+\rangle\subset \frak{h}^\ast$ das von $R^+$
erzeugte Untermonoid und $\leq$ die partielle Ordnung auf
$\frak{h}^\ast$ mit $\lambda\leq\mu\IFF \mu \in \lambda + |R^+\rangle.$ 
Jede kurze exakte Sequenz $\nabla_{T} (\mu) 
\hookrightarrow M \twoheadrightarrow \Delta_{T} (\lambda)$
mit $M \in \cal{O} (T)$ und $ \lambda \not\leq \mu$ 
spaltet, da jedes Urbild
in $M^{\lambda}$ des kanonischen Erzeugers von $\Delta_{T} (\lambda)$ bereits
von $\frak{n} $ annulliert werden mu"s und folglich eine Spaltung liefert.
Im Fall $\lambda \leq \mu$ gehen wir mit $d$ zur dualen
Situation "uber. Das zeigt die Trivialit"at der fraglichen 
Erweiterungen.
Den einfacheren Fall der R"aume von Homomorphismen behandelt man genauso.
\end{proof}

\begin{Korollar}\label{BCH}
Seien  $M,N \in \cal{O} (T).$ 
Ist $M$ ein direkter Summand eines Objekts mit endlicher  $\Delta_T$-Fahne und
$N$ ein direkter Summand eines Objekts mit endlicher $\nabla_T$-Fahne, so
ist
der Morphismenraum $\op{Hom}_{\cal{O}(T)} (M,N)$ ein endlich
erzeugter projektiver $T$-Modul und f"ur 
jede Erweiterung $T \ra T^{\prime}$ liefert
die offensichtliche Abbildung einen Isomorphismus
$$\op{Hom}_{\cal{O}(T)} (M,N) \otimes_{T} T^{\prime} 
\;\;\overset{\sim}{\ra}\;\; \op{Hom}_{\cal{O}(T^{\prime})}
(M \otimes_{T} T^{\prime}, N\otimes_{T} T^{\prime})$$
\end{Korollar}
\begin{proof}[Beweis]
Das folgt sofort aus \ref{EHo}. Die Details "uberlasse ich dem Leser.
\end{proof}

\begin{Bemerkung}\label{he}
Ist $Q\in\op{Kring}^S$ ein K"orper und landen f"ur alle Wurzeln $\alpha$ die 
Kowurzeln $\alpha^{\vee}$ unter
$S \ra Q$ nicht in $\Bbb{Z} \subset Q$, so ist die Kategorie 
$\cal{O}(Q)$ halbeinfach (d.h.\ alle Surjektionen
spalten) mit einfachen
Objekten $\Delta_{Q} (\lambda) = \nabla_{Q} (\lambda)$ f"ur 
$\lambda \in \frak{h}^{\ast}.$
\end{Bemerkung}
\section{Deformation unzerlegbarer Kippmoduln}
\begin{Bemerkung}
  Bezeichne $D=S_{(0)}$ der lokale Ring am Nullpunkt von $\frak{h}^{\ast}.$
  Ist $\lambda \in \frak{h}^{\ast}$ gegeben mit $\Delta (\lambda)$ einfach,
  so ist die kanonische Abbildung ein Isomorphismus
  $$\Delta_{D} (\lambda) \overset{\sim}{\ra} \nabla_{D} (\lambda)$$
  In der Tat
  reicht es zu zeigen, da"s diese Abbildung Isomorphismen auf allen
  Gewichtsr"aumen induziert, und diese sind freie Moduln von endlichem Rang
  "uber dem lokalen Ring $D.$ Es reicht also nach Nakayama zu zeigen,
  da"s unsere Abbildung unter $\otimes_{D} \Bbb{C}$ ein Isomorphismus wird, und
  das folgt sofort aus unserer Voraussetzung.
\end{Bemerkung}
\begin{Definition} Gegeben $T\in\op{Kring}^S$
bezeichne $\cal{K} (T) \subset \cal{O} (T)$ die kleinste 
Unterkategorie, die (1) diejenigen $\Delta_{T} (\lambda)$
umfa"st, f"ur die die kanonische Abbildung einen Isomorphismus 
$\Delta_{T} (\lambda) \overset{\sim}{\ra} \nabla_{T}
(\lambda)$ liefert, die (2) stabil ist unter dem Tensorieren mit 
endlichdimensionalen Darstellungen
von $\frak{g}$, die (3) stabil ist unter dem Bilden direkter 
Summanden und die (4) mit jedem Objekt auch alle
dazu isomorphen Objekte enth"alt.
Wir nennen $\cal{K} (T)$ die Kategorie der {\bf $T$-deformierten Kippmoduln}.
\end{Definition}
\begin{Bemerkung}
F"ur $\DC=\DC_0\in \op{Kring}^S$ sind
die Objekte von $\cal{K} (\DC)$ 
die "ublichen Kippmoduln
in der klassischen BGG-Kategorie $\cal{O}.$
\end{Bemerkung}
\begin{Proposition}\label{CTD}
  Ist $T\in\op{Kring}^S$ 
ein vollst"andiger lokaler Ring ``unter $S$''
derart, da"s das Urbild in $S$ seines maximalen Ideals gerade
das Verschwindungsideal des
  Nullpunkts von $\frak{h}^{\ast}$ ist, so induziert das Spezialisieren
  $$\otimes_{T} \Bbb{C} : \cal{K} (T)\ra\cal{K} (\Bbb{C})$$
  eine
  Bijektion auf Isomorphieklassen.
\end{Proposition}
\begin{proof}[Beweis]
Die Kippmoduln in $\cal{O}$ sind gerade die direkten Summanden von
Tensorprodukten einfacher Vermamoduln mit endlichdimensionalen
Darstellungen.
Die Proposition folgt mit \ref{BCH} 
aus den Definitionen und allgemeinen 
Tatsachen \cite{Ben}
"uber das Liften von Idempotenten.
\end{proof}
\begin{Bemerkung}
Gegeben $\lambda\in \frak{h}^\ast$ notieren wir
$K_T(\lambda)\in \cal{K} (T)$ die $T$-Deformation 
des unzerlegbaren Kippmoduls $K(\lambda)\in \cal{O}$ mit
h"ochstem Gewicht $\lambda.$
\end{Bemerkung}
\section{Die Andersen-Filtrierung}

\begin{Bemerkung}
Seien gegeben $K \in \frak{g}\op{-Mod_{\DC}-}T$ 
und
$\lambda \in \frak{h}^{\ast}.$ 
Der "Ubersichtlichkeit halber verwenden  wir 
hier die Abk"urzungen
$\Delta_{T} (\lambda)=\Delta,$ $\nabla_{T} (\lambda)=\nabla$
und $\op{Hom}_{\frak{g}-T}=\op{Hom}$ und
betrachten  die 
durch die Komposition gegebene $T$-bilineare Paarung 
$$\op{Hom} (\Delta, K) \times \op{Hom}
(K,\nabla) \ra \op{Hom} (\Delta,
\nabla) = T$$
Bezeichnen wir f"ur jeden $T$-Modul $H$ mit $H^\ast$ 
den $T$-Modul $\op{Hom}_T(M,T),$ so induziert unsere 
Paarung eine  Abbildung
$$E=E_\lambda(K):
\op{Hom} (\Delta , K) \ra \op{Hom}
(K, \nabla)^{\ast}$$
Ist hier $K$ ein Kippmodul, so geschieht unsere
Abbildung nach \ref{BCH} zwischen endlich erzeugten projektiven $T$-Moduln.
Ist zus"atzlich $T\in\op{Kring}^S$ ein Integrit"atsbereich und
  erf"ullt $Q = \op{Quot}T$ die Bedingung aus Bemerkung \ref{he}, ist also
  $\cal{O}(Q)$ halbeinfach mit einfachen Objekten $\Delta_{Q} (\lambda) =
  \nabla_{Q} (\lambda),$ so ist unsere Paarung nichtentartet "uber $Q$ und 
unsere Abbildung $E_\lambda(K)$
induziert einen Isomorphismus "uber $Q$ 
und ist insbesondere eine Injektion.  
  Ist speziell $T =\Bbb{C} [[t]]$ der Ring der formalen Potenzreihen l"angs
  einer Gerade $t \delta \subset \frak{h}^{\ast},$ die in keiner Spiegelebene
  der Weylgruppe enthalten ist, so erf"ullt $Q=\op{Quot}\Bbb{C} [[t]]$ 
die Bedingung aus
  Bemerkung \ref{he}. Ist dann $K\in \cal{K}(\Bbb{C} [[t]])$ 
ein deformierter Kippmodul, so
  k"onnen wir vermittels der Einbettung
  $$E_\lambda(K):\op{Hom} (\Delta , K) \hra
  \op{Hom} (K, \nabla)^{\ast}$$
  zwischen freien
  $\Bbb{C} [[t]]$-Moduln von endlichem Rang die 
offensichtliche Filtrierung auf der rechten
  Seite durch die $t^{i}\op{Hom} (K, \nabla)^{\ast}$ nach links
zur"uckziehen und so
  eine Filtrierung auf $\op{Hom} (\Delta, K)
=\op{Hom}_{\frak{g}-\Bbb{C} [[t]]} (\Delta_{\Bbb{C} [[t]]} (\lambda), K)$
  erhalten.
\end{Bemerkung}

\begin{Definition}\label{DAF}
Ist $K_\DC\in\cal{K}(\DC)$ ein Kippmodul der "ublichen Kategorie
$\cal{O}$ und $K\in \cal{K}(\DC[[t]])$ eine $\DC[[t]]$-Deformation 
im Sinne von \ref{CTD} mit $S\ra\DC[[t]]$ 
der Restriktion auf eine formale Umgebung des
Nullpunkts in der Gerade $t\rho$ mit $\rho$ wie in
\ref{rhod},
so nennen wir das Bild der eben erkl"arten 
Filtrierung unter der Spezialisierung
  $\otimes_{\Bbb{C}[[t]]} \Bbb{C}$ 
die {\bf Andersen-Filtrierung auf}
  $\op{Hom}_{\frak{g}} (\Delta (\lambda), K_{\Bbb{C}}).$
\end{Definition}
\begin{Bemerkung}
Es bleibe dem Leser "uberlassen, die Unabh"angigkeit dieser 
Filtrierung von der Wahl der nur bis auf nichteindeutigen
Isomorphismus wohldefinierten Deformation unseres Kippmoduls 
nachzuweisen. 
Das Ziel dieser Arbeit ist die Berechnung der Dimensionen
der Subquotienten der Andersen-Filtrierung auf
$\op{Hom}_{\frak{g}} (\Delta (\lambda), K(\mu))$
f"ur alle $\lambda,\mu\in\frak{h}^\ast$ oder vielmehr ihre
Beschreibung durch Koeffizienten von Kazhdan-Lusztig-Polynomen.
Die endg"ultige Formel f"ur die
Dimension des $i$-ten Subquotienten einer
beliebigen Andersen-Filtrierung lautet
\begin{displaymath}
\dim_\DC \bar{F}^i\op{Hom}_{\frak{g}}(\Delta (\lambda_{\bar{y}}), 
K(\lambda_{\bar{x}}))=h_{y,x}^i
\end{displaymath}
Hierbei gehen wir aus von einem
$\rho$-dominanten Gewicht $\lambda\in\frak{h}^\ast_{\op{dom}}$ im
Sinne von \ref{rhod}. Es liefert zwei Untergruppen
$W_{\bar{\lambda}}\supset W_\lambda$ der Weylgruppe wie in
\ref{WLL} ausgef"uhrt wird, und $\bar{x},\bar{y}$ meinen Nebenklassen
aus $W_{\bar{\lambda}}/ W_\lambda$ mit $x,y$ als l"angsten 
Repr"asentanten. Schlie"slich meint 
$\lambda_{\bar{x}}=w_{\bar{\lambda}}\bar{x}\cdot\lambda$ wie in
\ref{NZ} mit $w_{\bar{\lambda}}$ dem l"angsten Element von
$W_{\bar{\lambda}},$ und $h^i_{y,x}$ wird dadurch charakterisiert,
da"s in den Notationen von Kazhdan und Lusztig das KL-Polynom
$P_{y,x}$ in Bezug auf die Coxetergruppe $W_{\bar{\lambda}}$ mit
ihrer L"angenfunktion $l$ gegeben wird durch die Formel
$$P_{y,x}=\sum h^i_{y,x} q^{(l(x)-l(y)-i)/2}$$
In Wirklichkeit zeigen wir, da"s die fragliche Filtrierung 
mit der Graduierungsfiltrierung zusammenf"allt, die von der
graduierten Version der Kategorie $\cal{O}$ induziert wird, 
aber es schien mir nicht m"oglich, diese Erkenntnis im Rahmen
eines Zeitschriftenartikels verst"andlich darzustellen.
\end{Bemerkung}
\begin{Bemerkung}
Die Jantzen-Filtrierung auf einem Vermamodul 
$\Delta (\lambda)$ induziert
nat"urlich Filtrierungen auf
$\op{Hom}_{\frak{g}} (P(\mu), \Delta (\lambda))$ f"ur
$P (\mu) \twoheadrightarrow \Delta (\mu)$
die projektive Decke in $\cal{O}.$
Diese Filtrierungen hinwiederum kommen in derselben Weise her von 
den durch $\Delta_{\Bbb{C} [[t]]}(\lambda)\ra\nabla_{\Bbb{C} [[t]]}(\lambda)$ 
induzierten Einbettungen
$$\op{Hom}_{\frak{g}-\Bbb{C} [[t]]}(P_{\Bbb{C} [[t]]} (\mu), 
\Delta_{\Bbb{C} [[t]]} (\lambda)) \hookrightarrow
\op{Hom}_{\frak{g}-\Bbb{C} [[t]]} (P_{\Bbb{C} [[t]]}(\mu), 
\nabla_{\Bbb{C} [[t]]} (\lambda))$$
mit 
$P_{\Bbb{C} [[t]]}(\mu)\sra \Delta_{\Bbb{C} [[t]]} (\lambda)$ 
den projektiven Decken 
in $\cal{O}(\Bbb{C} [[t]]).$ Das zeigt die Analogie 
zwischen beiden Filtrierungen. Ich erwarte, da"s Archipov's Kippfunktor
sogar diese Filtrierungen identifiziert, kann es aber leider nicht beweisen,
ohne die Jantzen-Vermutung vorauszusetzen.
\end{Bemerkung}
\section{Deformierte Verschiebung}
\begin{Bemerkung}
  Bezeichnet $Z \subset U(\frak{g})$ das Zentrum, so operiert nat"urlich $Z
  \otimes T$ auf jedem Bimodul $M \in \frak{g}\op{-Mod_{\DC}-}T $. 
Wir betrachten nun das push-out-Diagramm von $\Bbb{C}$-Algebren
  $$\begin{array}{ccccc}
    & & Z \otimes T & &\\
    &\nearrow & & \searrow &\\
    Z \otimes \Bbb{C} [\frak{h}^{\ast}] & & & & \Bbb{C}[\frak{h}^{\ast}]
    \otimes T\\
    & \searrow & \ & \nearrow\\
    & &\Bbb{C}[\frak{h}^{\ast}] \otimes \Bbb{C} [\frak{h}^{\ast}]&&
 \end{array}$$
 wobei wir f"ur $\xi : Z \ra \Bbb{C} [\frak{h}^{\ast}]$ stets die Variante des
 Harish-Chandra-Homo\-mor\-phismus mit $\xi (z) -z \in U \frak{n} $ nehmen. 
Er liefert eine endliche Ringerweiterung und dasselbe gilt dann auch 
f"ur beide nach unten gerichteten Pfeile unseres Diagramms.
 Der Graph der Addition mit einem $\lambda \in \frak{h}^{\ast}$ ist eine
irreduzible 
abgeschlossene Teilmenge von $\op{Spec} (\Bbb{C}[\frak{h}^{\ast}] \otimes
 \Bbb{C} [\frak{h}^{\ast}])$ und dasselbe gilt f"ur sein Bild in
$\op{Spec} (Z \otimes \Bbb{C} [\frak{h}^{\ast}]).$ 
Das Urbild in $\op{Spec} (Z\otimes T)$ dieses
 Bildes bezeichnen wir mit
 $\Xi_{\lambda}\As \op{Spec} (Z \otimes T)$.  
Per definitionem haben $\Delta_{T}
 (\lambda)$ und $\nabla_{T} (\lambda)$ beide als $Z\otimes T$-Moduln Tr"ager in
 $\Xi_{\lambda}$.
\end{Bemerkung}
\begin{Lemma}
Der Tr"ager in
$\op{Spec} (Z \otimes T)$
jedes Elements eines Moduls 
$M \in \cal{O} (T)$ ist enthalten in einer endlichen
Vereinigung von Mengen der Gestalt $\Xi_{\lambda}$ 
mit $ \lambda \in \frak{h}^{\ast}.$
\end{Lemma}
\begin{proof}[Beweis]
Sei $v$ unser Element.
Wir d"urfen annehmen, da"s gilt $v \in M^{\lambda}$ 
f"ur ein $\lambda \in \frak{h}^{\ast}.$
Wir d"urfen weiter annehmen,
da"s das Erzeugnis von $v$ enthalten ist in 
$\bigoplus_{\mu \leq \nu}M^{\lambda + \mu}$
f"ur ein dominantes ganzes Gewicht $\nu \in X^{+}.$
Dann besitzt 
$$U (\frak{g}) \otimes_{U(\frak{b})} \tau_{\leq \lambda +\nu} \left(U(\frak{b})
\otimes_{U(\frak{h})} (\Bbb{C}_{\lambda} \otimes T)\right)$$
mit hoffentlich selbsterkl"arendem $\tau_{\leq \lambda + \nu}$ 
eine endliche $\Delta_{T}$-Fahne
und unser $v$ liegt im Bild eines Homomorphismus von 
besagtem Objekt nach $M$.
\end{proof}

\begin{Definition}\label{rhod}
Bezeichne $\rho=\rho(R^+)$ die Halbsumme der positiven Wurzeln.
Wir setzen
$$\frak{h}^\ast_{\op{dom}} = \{ \lambda \in \frak{h}^\ast \mid
\langle \lambda + \rho, \alpha^{\vee} \rangle \not\in \{ -1,-2,
\ldots\} \; \forall \alpha \in R^{+}\}$$
und nennen
die  Elemente dieser Menge die \defind{$\rho$-dominanten Gewichte}.
Wir verwenden die "ubliche Notation $w\cdot\lambda=w(\lambda+\rho)-\rho$
f"ur die zum Fixpunkt $-\rho$ verschobene Operation der Weylgruppe.
\end{Definition}

\begin{Satz}[Zerlegung deformierter Kategorien]
Sei $T$ ein $S_{(0)}$-Ring, d.h.\ der
Morphismus $S\ra T$ faktorisiere "uber den lokalen Ring $S_{(0)}$
von $\frak{h}^\ast$ am Ursprung.
So haben wir eine Zerlegung
$$\cal{O} (T) = \prod_{\lambda \in \frak{h}^{\ast}_{\op{dom}} }
\cal{O}_{\lambda}(T)$$
wobei $\cal{O}_{\lambda} (T)$ aus allen Objekten $M \in \cal{O}(T)$ besteht mit
$M = \bigoplus_{\nu \in \lambda + \Bbb{Z} R} M^{\nu} $ und $
\op{supp}_{Z \otimes T} M
\subset \bigcup_{w \in W} \Xi_{w \cdot \lambda}.$
\end{Satz}
\begin{proof}[Beweis]
Aus $\Xi_{\lambda} \cap \Xi_{\mu} \neq \emptyset$ folgt 
f"ur $T=S_{(0)}=D$ bereits $W \cdot 
\lambda = W \cdot \mu.$
Der Rest des Arguments l"auft wie im 
Fall $T = \Bbb{C}.$
\end{proof}
\begin{Bemerkung}
Ebenso wie im nichtdeformierten Fall haben wir f"ur $\lambda, 
\mu \in \frak{h}_{\op{dom}}^{\ast}$
mit ganzer Differenz $\lambda -\mu \in X$ Verschiebungsfunktoren
$$T^{\mu}_{\lambda} : \cal{O}_{\lambda} (T) \ra \cal{O}_{\mu} (T)$$
die exakt sind und Adjunktionen $(T^{\mu}_{\lambda}, T^{\lambda}_{\mu})$ 
erf"ullen und
auch sonst die "ublichen Eigenschaften haben.
Wir nennen sie die 
{\bf deformierten Verschiebungen}.
Die Kategorie der deformierten Kippmoduln in einem unserer Bl"ocke
notieren wir $\cal{K}T)\cap \cal{O}_\lambda(T)=\cal{K}_\lambda(T).$
\end{Bemerkung}
\begin{Bemerkung}\label{WLL}
Ist $T$ eine $D$-Algebra, so
ist f"ur $\lambda \in \frak{h}^{\ast}_{\op{dom}}$ 
der deformierte Verma-Modul $\Delta_{T} (\lambda)$
projektiv in $\cal{O}_{\lambda} (T)$. Die Isotropiegruppe eines Gewichts 
$\lambda\in\frak{h}^{\ast}$  unter
der dot-Operation der Weylgruppe notieren wir $W_\lambda,$ 
die Isotropiegruppe seiner Nebenklasse $\bar{\lambda}=\lambda+\langle R\rangle$
unter dem Wurzelgitter $W_{\bar{\lambda}}.$ 
Das l"angste Element von  $W_{\bar{\lambda}}$ notieren wir 
$w_{\bar{\lambda}},$ 
die Ringe der Invarianten f"ur die  
nat"urliche Operation der Gruppen $W_\lambda\subset W_{\bar{\lambda}}$ auf $D$ 
bezeichnen wir mit $D^\lambda\supset D^{\bar{\lambda}}.$
\end{Bemerkung}
\begin{Satz}[Deformation projektiver Objekte]
Der Funktor $\otimes_{D} \Bbb{C} : 
\cal{O} (D) \ra \cal{O} (\Bbb{C})$ induziert
eine Bijektion zwischen den Isomorphieklassen endlich erzeugter 
projektiver Objekte in beiden
Kategorien.
\end{Satz}
\begin{proof}[Beweis]
\cite{So-A}.
\end{proof}
\begin{Definition}
Gegeben $\lambda \in \frak{h}^{\ast}$ bezeichne 
$P_{D} (\lambda) \in \cal{O} (D)$ das
endlich erzeugte projektive Objekt, das unter $\otimes_{D} \Bbb{C} $ zu
$P (\lambda)$ spezialisiert. Wir nennen es die {\bf Deformation des
Projektiven $P (\lambda)$.} Gegeben $\lambda \in \frak{h}^{\ast}_{\op{dom}}$ 
benutzen wir f"ur den deformierten antidominanten Projektiven die Abk"urzung
$P_{D} (w_{\bar{\lambda}} \cdot\lambda)=A_D(\lambda)=A(\lambda).$ 
\end{Definition}
\begin{Satz}[Endomorphismen antidominanter Projektiver]
Gegeben $\lambda \in \frak{h}^{\ast}_{\op{dom}}$ induziert 
die Multiplikation eine
Surjektion
$Z \otimes D \twoheadrightarrow 
\op{End}_{\cal{O}(D)} A (\lambda).$
Bezeichnet $(+\lambda)^{\sharp} : S \ra S$ den 
Komorphismus zu $(+\lambda) : \frak{h}^{\ast}\ra
\frak{h}^{\ast}$, so hat die Komposition
$$Z \otimes D \overset{\xi \otimes \op{id}}{\lra} 
S \otimes D
\overset{(+\lambda)^{\sharp}\otimes\op{id}}{\lra} S 
\otimes D \ra D \otimes_{D^{\bar{\lambda}}}
D$$
das Bild $D^{{\lambda}} \otimes_{D^{\bar{\lambda}}} D$ und denselben
Kern wie die Surjektion aus dem ersten Teil und wir erhalten so 
einen Isomorphismus
$$D^{{\lambda}} \otimes_{D^{\bar{\lambda}}} D \overset{\sim}{\ra}
\op{End}_{\cal{O}(D)} A ( \lambda)$$
\end{Satz}
\begin{proof}[Beweis]
F"ur $\lambda$ ganz wird das in \cite{HCH} gezeigt.
Der Beweis im Allgemeinen ist im Wesentlichen derselbe.
\end{proof}
\begin{Bemerkung}
Der "Ubersichtlichkeit halber verwenden wir im Folgenden
meist 
die Abk"urzungen 
$\op{Hom}_{\cal{O}(D)}=\op{Hom}$ und $\op{End}_{\cal{O}(D)}=\op{End}.$
Jede Wahl eines deformierten antidominanten
Projektiven  $A ( \lambda)$ f"ur
$\lambda \in \frak{h}^{\ast}_{\op{dom}}$ liefert  
vermittels der Vorschrift
$\Bbb{V} = \Bbb{V}_{D} 
= \op{Hom}_{\cal{O}(D)}(A( \lambda),\;)$   
einen exakten
Funktor
$$
\Bbb{V}:\cal{O} (D) \ra D^{{\lambda}}\op{-Mod_{\DC}-}D$$
der nat"urlich nur auf $\cal{O}_{\lambda} (D)$ von Null verschieden ist.
Ist weiter $\mu \in \frak{h}^{\ast}_{\op{dom}}$ gegeben mit
$\lambda -\mu \in X$ und $W_{\mu} \supset W_{\lambda}$ und w"ahlen
wir einen Isomorphismus $T^{\lambda}_{\mu} A( \mu) \overset{\sim}{\ra}
A(\lambda)$, so kommutiert das Diagramm
$$\begin{array}{ccc}
D^{{\mu}} \otimes_{D^{\bar{\mu}}} D& 
\overset{\sim}{\ra} & \op{End}
A(\mu)\\
\downarrow & & \downarrow \\
D^{{\lambda}} \otimes_{D^{\bar{\mu}}} D & 
\overset{\sim}{\ra} & \op{End}
A(\lambda)
\end{array}$$
mit der von der Einbettung $D^{{\mu}} \subset D^{{\lambda}}$ 
induzierten linken
Vertikale und der durch $T^{\lambda}_{\mu}$ und unseren Isomorphismus 
induzierten rechten
Vertikale, siehe \cite{HCH}.  Halten wir so einen Isomorphismus fest und
  w"ahlen au"serdem eine Adjunktion $(T^{\lambda}_{\mu}, T^{\mu}_{\lambda})$, 
so
  erhalten wir  Isomorphismen
  $$\op{Hom} (A(\mu), T^{\mu}_{\lambda} M)
  \overset{\sim}{\ra} 
\op{Hom} (T^{\lambda}_{\mu} A (
  \mu), M) \overset{\sim}{\ra} 
\op{Hom} (A( \lambda),M)$$
  die zusammen eine "Aquivalenz von Funktoren definieren, bis auf die das
  Diagramm
  $$\begin{array}{ccc}
    \cal{O}_{\lambda} (D) & \overset{\Bbb{V}}{\ra} & D^{{\lambda}} 
\op{-Mod_{\DC}-}D\\[2mm]
    T_{\lambda}^{\mu} \downarrow & & \downarrow \op{res}\\[2mm]
    \cal{O}_{\mu} (D) & \overset{\Bbb{V}}{\ra} & D^{{\mu}}\op{-Mod_{\DC}-}D
\end{array}$$
kommutiert.  Mithilfe der Adjunktionen finden wir auch eine "Aquivalenz von
Funktoren, bis auf die das Diagramm
$$\begin{array}{ccl}
  \cal{O}_{\mu} (D) 
& \overset{\Bbb{V}}{\ra} & D^{{\mu}} \op{-Mod_{\DC}-}D\\[2mm]
  T^{\lambda}_{\mu} \downarrow & 
& \;\;\;\;\;\;\;\;\downarrow D^{{\lambda}}
  \otimes_{D^{{\mu}}}\\[2mm]
  \cal{O}_{\lambda}(D) &\overset{\Bbb{V}}{\ra} & D^{{\lambda}}
\op{-Mod_{\DC}-}
  D
\end{array}$$
kommutiert.  
\end{Bemerkung}

\begin{Bemerkung}
Gegeben $\lambda, \mu \in \frak{h}^*_{\op{dom}}$
mit ganzer Differenz kann man allgemeiner die Verschiebungen
$
T^{\mu}_{\lambda} : \mathcal{O}_{\lambda} (D) \ra \mathcal{O}_{\mu} (D)$
betrachten, die sich bekanntlich auch schreiben lassen als $T^{\mu}_{\lambda}
\cong T_{\nu}^{\mu} T^{\nu}_{\lambda}$ f"ur ein und 
jedes $\nu \in \frak{h}^*_{\op{dom}}$
mit ganzer Differenz zu $\lambda$ und $\mu$ und der 
Eigenschaft $W_{\nu} = W_{\lambda}
\cap W_{\mu}$.
Bilden wir $D^{\mu}_{\lambda} = D^{W_{\lambda} \cap W_{\mu}} \in D^{\mu}
\op{-Mod_{\Bbb{C}}-}D^{\lambda}$, so k"onnen wir unsere 
Diagramme zusammenfassen
zum bis auf nat"urliche "Aquivalenz kommutierenden funktoriellen Diagramm
\begin{displaymath}
\xymatrix{
\mathcal{O}_{\lambda}(D) \ar[r]^-{\Bbb{V}} \ar[d]_{T^{\mu}_{\lambda}} 
& D^{\lambda}
\op{-Mod_{\Bbb{C}}-}D \ar[d]^{D^{\mu}_{\lambda} \otimes_{D^{\lambda}}}\\
\mathcal{O}_{\mu} (D) \ar[r]^-{\Bbb{V}} & D^{\mu}\op{-Mod_{\Bbb{C}}-}D
}
\end{displaymath}
Gehen wir zu den Adjungierten der Vertikalen "uber, so erhalten wir auch ein
bis auf nat"urliche "Aquivalenz von Funktoren kommutierendes Diagramm
\begin{displaymath}
\xymatrix{
\mathcal{O}_{\mu} (D) \ar[r]^-{\Bbb{V}} \ar[d]_{T^{\lambda}_{\mu}} &
D^{\mu}\op{-Mod_{\Bbb{C}}-}D\ar[d]^{\op{Hom}_{D^{\mu}} (D_{\lambda}^{\mu},
\quad) }\\
\mathcal{O}_{\lambda} (D) \ar[r]^-{\Bbb{V}} & D^{\lambda}\op{-Mod_{\Bbb{C}}-}
D
}
\end{displaymath}

\end{Bemerkung}

\begin{Theorem}[Struktursatz f"ur deformierte Kippmoduln]
Die Funktoren $\Bbb{V}$ sind volltreu auf deformierten Kippmoduln. 
Genauer gilt f"ur beliebiges $\lambda\in\frak{h}^\ast_{\op{dom}}$:
\begin{enumerate}
\item
Gegeben $K \in \cal{K}_\lambda (D) $ 
und $F \in \cal{O}_{\lambda} (D)$
ein Objekt mit  $\Delta_D$-Fahne induziert $\Bbb{V}$ einen Isomorphismus
$$\op{Hom}_{\frak{g}-D} (F, K) \overset{\sim}{\ra} 
\op{Hom}_{D^{\lambda}-D} (\Bbb{V} F, \Bbb{V}K)$$
\item
Gegeben $K \in \cal{K}_\lambda (D) $ 
und $F \in \cal{O}_{\lambda} (D)$ ein Objekt mit  $\nabla_D$-Fahne 
induziert
$\Bbb{V}$ einen Isomorphismus
$$\op{Hom}_{\frak{g}-D} (K,F) \overset{\sim}{\ra} \op{Hom}_{D^{\lambda}-D}
(\Bbb{V} K, \Bbb{V} F)$$
\end{enumerate}
\end{Theorem}
\begin{Bemerkung}
In gr"o"serer Allgemeinheit wird die erste Aussage als
Theorem 10 in \cite{FieC} gezeigt: Die Funktoren  $\Bbb{V}$ sind 
sogar volltreu auf beliebigen Objekten mit einer 
endlichen Filtrierung, bei der alle Subquotienten deformierte
Vermamoduln sind.
\end{Bemerkung}
\begin{Bemerkung}
Im nichtdeformierten Fall $T = \Bbb{C}$ ist der Funktor $\Bbb{V}$
volltreu auf den Kippmoduln eines gegebenen Blocks. In der Tat liefert ja f"ur
jedes maximale Ideal $\chi \subset Z$ und beliebig
vorgegebene  projektive Funktoren
$F, G : U / \chi U \op{-mod} \ra U \op{-mod}$ und einen beliebigen
Vermamodul $\Delta$ mit $\chi \Delta =0$
das Anwenden auf $\Delta$ eine Bijektion
\begin{displaymath}
\op{Trans}_{U/\chi U \op{-mod}\ra} (F,G) 
\overset{\sim}{\rightarrow} \op{Hom}_{\frak{g}}
(F\Delta, G\Delta)
\end{displaymath}
Das wird f"ur projektive Vermamoduln 
bewiesen in \cite{BG}, und da nach \cite{BGG-Aff}
die Einh"ullende surjektiv auf die 
$\op{ad}$-endlichen Endomorphismen jedes Vermamoduls
geht, funktioniert der dort gegebene 
Beweis allgemeiner f"ur jeden Vermamodul.
Die Einbettung eines einfachen Vermamoduls $\Delta_e$ in einen projektiven
Vermamodul $\Delta_p$ liefert also Bijektionen 
\begin{displaymath}
\op{Hom}_{\frak{g}} (F \Delta_e, G \Delta_e) 
\overset{\sim}{\rightarrow} \op{Hom}_{\frak{g}}
(F\Delta_p, G\Delta_p)
\end{displaymath}
Da sie auch Bijektionen 
$\Bbb{V} F \Delta_e \overset{\sim}{\rightarrow} \Bbb{V}
F \Delta_p$ liefert, folgt die Behauptung.
Im nichtdeformierten Fall kann jedoch von Volltreuheit 
auf Morphismen von Kippmoduln
zu dualen Vermamoduln oder von Vermamoduln zu Kippmoduln keine Rede sein.
\end{Bemerkung}

\begin{proof}[Beweis]
Gegeben $\lambda, \mu \in \frak{h}^{\ast}_{\op{dom}}$ mit
ganzer Differenz sei
$D^{\lambda}_{\mu} \in D^{\lambda} \op{-Mod_{\DC}-} D^{\mu}$ der Bimodul
$D^{W_{\lambda} \cap W_{\mu}}.$ Die vorhergehenden 
"Uberlegungen zeigen, da"s das Diagramm
$$\begin{array}{lcr}
\op{Hom}_{\frak{g}-D} (F,T^{\mu}_{\lambda}K)& \ra
&\op{Hom}_{\frak{g}-D} (T_{\mu}^{\lambda} F, K)\\
\;\;\;\;\;\;\;\;\;\;\;\;\;\;\;\;\;\;\da&&\da\;\;\;\;\;\;\;\;\;\;\;
\;\;\;\;\;\;\;\\
\op{Hom}_{D^{\mu}-D} (\Bbb{V}F, \Bbb{V}
T^{\mu}_{\lambda} K) & 
&\op{Hom}_{D^{\lambda}-D} (\Bbb{V} T^{\lambda}_{\mu}
F, \Bbb{V} K)\\
\;\;\;\;\;\;\;\;\;\;\;\;\;\;\;\;\;\;\da&&\da\;\;\;\;\;\;\;\;\;
\;\;\;\;\;\;\;\;\;\\
 \op{Hom}_{D^{\mu}-D} 
(\Bbb{V}F, \op{Hom}_{D^{\lambda}}(D^{\lambda}_{\mu},
\Bbb{V} K))& \ra
&\op{Hom}_{D^{\lambda}-D} 
(D^{\lambda}_{\mu} \otimes_{D^{\mu}}
\Bbb{V} F, \Bbb{V} K)
\end{array}$$
kommutiert, wenn
wir die beiden unteren vertikalen Morphismen 
mithilfe der eben eingef"uhrten  nat"urlichen "Aquivalenzen
erkl"aren und die waagerechten Morphismen mithilfe der Adjunktionen.
In diesem Diagramm sind alle Morphismen mit Ausnahme der beiden oberen
Vertikalen offensichtlich Isomorphismen. 
Ist die rechte obere Vertikale ein Isomorphismus, so mithin auch die 
linke obere Vertikale.
Gilt in anderen Worten 
unsere Behauptung f"ur $K$, so auch f"ur $T_{\lambda}^{\mu} K$.
Damit m"ussen wir sie nur f"ur $K$ einen deformierten einfachen Vermamodul
pr"ufen. Indem wir die Vermafahne von oben abarbeiten, d"urfen wir sogar
annehmen, da"s $F$ eine direkte Summe von Kopien eben dieses einfachen
Verma-Moduls ist.
In dem Fall ist aber die erste Behauptung klar. Die zweite Behauptung zeigt man
analog.
\end{proof}

\section{Geometrische Argumente}

\begin{Notation}
  Bezeichne $\op{gMod-}\! A$ die Kategorie der graduierten Rechtsmoduln "uber
  einem graduierten Ring $A.$ Bezeichne $\op{Der}_G(X)$ bzw.\ 
$\op{Der}_G^+(X)$ die "aquivariante bzw.\ die "aquivariante 
  gegen die Pfeile beschr"ankte derivierte Kategorie zu einer komplexen
  algebraischen Variet"at $X$ mit einer Operation einer komplexen algebraischen
  Gruppe $G$ und bezeichne $\op{Der}_G(\cal{F},\cal{G})$ die Morphismen in
  diesen Kategorien.
\end{Notation}

\begin{Bemerkung}\label{AGSU}
 Wir arbeiten im folgenden in der Kohomologie stets mit 
  komplexen Koeffizienten.  
Seien ganz allgemein  
$X$ eine komplexe algebraische Variet"at mit einer Operation
einer algebraischen Gruppe $B.$ 
Sei $X=\coprod_{a\in A}X_a$ eine Stratifizierung in 
irreduzible lokal abgeschlossene glatte 
$B$-stabile Untervariet"aten derart, da"s der 
Abschlu"s jedes Stratums eine Vereinigung von Strata ist. 
Bezeichne $|a|$ die Dimension von $X_a$ und
$\cal{C}_a=\underline{X_a}[|a|]$ 
die ``konstante perverse Garbe'' in $\op{Der}_B(X_a).$
Bezeichne weiter $j_a:X_a\hra X$ die Einbettung.
Seien nun 
$\mathcal{F}, \mathcal{G} \in
\op{Der}_{B} (X)$ gegeben mit der Eigenschaft, da"s f"ur alle
$a\in A$ gilt
\begin{displaymath}
\begin{array}{ccl}
j^{\ast}_{a} \mathcal{F} & \cong & 
\bigoplus_{\nu} f^\nu_{a} \mathcal{C}_a [\nu]\\[2mm]
j^!_{a} \mathcal{G} & \cong & 
\bigoplus_{\nu} g^{\nu}_{a} \mathcal{C}_a [\nu]
\end{array}
\end{displaymath}
in $\op{Der}_{B} (X_a)$ f"ur geeignete 
$f^{\nu}_a,g^{\nu}_a \in \Bbb{N}.$
Haben wir unter diesen Annahmen zus"atzlich 
$f^\nu_{a}= 0=g^\nu_{a}$ f"ur $\nu+|a|$ 
ungerade und $H^\nu_B(X_a)= 0$ f"ur $\nu$ ungerade, 
so induziert $\Bbb{H}_B$ f''ur alle $\ast$ eine Injektion
$$\op{Der}_B(\cal{F},\cal{G}[\ast])\hra 
\op{Hom}(\Bbb{H}_B\cal{F}, \Bbb{H}_B\cal{G})$$
und die Dimensionen der homogenen Komponenten auf  
der linken Seite werden gegeben durch die Formel
$$\dim\op{Der}_B(\cal{F},\cal{G}[n])=
\sum_{\nu-\mu+k=n,\; a\in A} f^\nu_a g^\mu_a \dim H^k_B(X_a)$$
Der Beweis l"auft v"ollig analog zum Beweis von 
Proposition 3 auf Seite 404 von \cite{So-L}
und soll hier nicht wiederholt werden.
\end{Bemerkung}

\begin{Bemerkung}
  Seien $G \supset P=P_\iota\supset B \supset T$ 
eine halbeinfache komplexe algebraische Gruppe, eine Parabolische,
  eine Borel und ein maximaler Torus. 
Seien  $W_\iota\subset W$ die Weylgruppen von $P\subset G$ und sei
$L\supset T$ die Levi von $P.$
Wir lassen $B\times P$ operieren auf $G$ vermittels
der Vorschrift $(b,p)g=bgp^{-1}.$
Bekanntlich wird der "aquivariante Kohomologiering
  $H^{\ast}_{B\times P}(G)$ mit dem Zur"uckholen ein Quotient von
  $H^\ast_{B \times P} (\op{pt}) = H^\ast_{T\times L} (\op{pt}) =
  R\otimes_{\Bbb{C}} R^\iota$ f"ur
$R=\cal{O}(\op{Lie}T)$ der Ring der regul"aren Funktionen auf
  $\op{Lie}T$, graduiert durch die Bedingung, da"s 
Linearformen 
  homogen vom Grad Zwei sein sollen,
und $R^\iota$ die $W_\iota$-Invarianten in $R.$ 
Wir erhalten damit einen kanonischen
  Isomorphismus
\begin{displaymath}
c : R \otimes_{R^W} R^\iota\overset{\sim}{\rightarrow} 
H^\ast_{B\times P} (G)
\end{displaymath}
So
liefert die "aquivariante Hyperkohomologie
\begin{displaymath}
\Bbb{H}^\ast_{B \times P} : \op{Der}_{B \times P}^+ (G) 
\ra \op{gMod-} H^\ast_{B \times P} (\op{pt})
\end{displaymath}
unter unserer Identifikation des "aquivarianten Kohomologierings und der
Identifikation 
$\op{Der}_{B\times P}^+ (G)\cong \op{Der}_B^+(G/P)$ einen Funktor in
die $\Bbb{Z}$-graduierten $R$-$R^\iota$-Bimoduln
\begin{displaymath}
\Bbb{H}_B=\Bbb{H}_B^\ast: \op{Der}^+_{B} (G/P) \ra R\op{-gMod-} R^\iota
\end{displaymath}  
Betrachten wir in  $\op{Der}_{B}^{+} (G/P)$ f"ur $x \in W/W_\iota$ 
nun den Schnittkohomologiekomplex
  $\cal{IC}_x$ zum Abschlu"s von ${B x P/P}.$ Sei 
$\cal{C}_y $ die konstante perverse 
Garbe auf $B
  y B/P,$ die also als Komplex von gew"ohnlichen Garben im Grad $-l(y)$
  konzentriert ist, und bezeichne $j_y: By B/P
  \hookrightarrow G/P$ die Einbettung.
\end{Bemerkung}
\begin{Satz}\label{VTM}
Der Funktor $\Bbb{H}_B$ ist volltreu f"ur Morphismen $\cal{IC}_x \ra
j_{y\ast}\cal{C}_y [n]$ und  $j_{x!}\cal{C}_y  \ra \cal{IC}_y[n]$
und $\cal{IC}_x \ra\cal{IC}_y[n]$
in $\op{Der}_{B}^+ (G/P).$
\end{Satz}
\begin{proof}
In \cite{So-L}, Proposition 2, Seite 402
wird der Satz f"ur Homomorphismen
$\cal{IC}_x \ra\cal{IC}_y[n]$ und $P=B$ bewiesen.
Ich werde im folgenden darlegen, in welcher Weise der 
dort gegebene Beweis
mit eher unwesentlichen "Anderungen auch 
diesen allgemeineren Satz zeigt.
Zun"achst beschr"anken wir uns auf den Fall $P=B.$
Nach Lemma 6 auf Seite 405 von loc.cit.\ 
in Verbindung mit \ref{AGSU} ist der Funktor im Lemma schon mal treu
und die Dimension der fraglichen Hom-R"aume ist bekannt. 
Nach Resultaten in \cite{SoBi} kennen wir jedoch auch 
die Dimensionen der Hom-R"aume im Bild und damit k"onnen wir 
das Argument mit einem Dimensionsvergleich abschlie"sen.
Genauer ergibt sich mit \ref{AGSU}   die Formel
\begin{displaymath}
\op{dim}_{\Bbb{C}} 
\op{Der}_{B} (\mathcal{I}\mathcal{C}_x, j_{y\ast} \mathcal{C}_y [n])
=\sum_{k+i = n} n^i_{y,x} \op{dim}_{\Bbb{C}} H^k_{B} (ByB/B)
\end{displaymath}
Hier werden die $n^i_{y,x}$ gegeben als 
Koeffizienten von Kazhdan-Lusztig-Polynomen
und es gilt genauer
\begin{displaymath}
\sum_{y,i} n^i_{y,x} q^{-i/2} \tilde{T}_{y} = C^{\prime}_{x}
\end{displaymath}
in Lusztig's Notationen f"ur die Hecke-Algebra alias
$\sum_{y,i} n^i_{y,x} v^i H_y = \underline{H}_x$ in den 
Notationen aus \cite{So-K}.
Andererseits wird in loc.cit.\ 
Lemma 5, Seite 402
f"ur $P=B$  bewiesen, da"s die $\Bbb{H}_B \cal{IC}_{x}$ 
gerade die speziellen Bimoduln
  $$\Bbb{H}_B \cal{IC}_{x}\cong B_x$$ 
sind, die in \cite{HCH} und \cite{SoBi} diskutiert werden.  
Nun  erinnern wir
 an den graduierten Bimodul $R_y$ aus \cite{SoBi}, der frei ist
vom Rang Eins
   von rechts und von links mit ein- und demselben 
Erzeuger $1_y$ im Grad Null und der
  Eigenschaft $r1_y=1_y r^y$ f"ur $r^y=y^{-1}(r),$ sowie an seine beiden in der
  Graduierung verschobenen Versionen $\Delta_y=R_y[-l(y)]$ und
  $\nabla_y=R_y[l(y)].$ 
Man zeigt unschwer 
$\Bbb{H}_Bj_{y\ast}\cal{C}_y\cong \nabla_y$ und 
$\Bbb{H}_Bj_{y!}\cal{C}_y\cong \Delta_y.$
Wir m"ussen also die Gleichheit
der Dimensionen
\begin{displaymath}
\op{dim}_{\Bbb{C}} \op{Der}_{B} 
(\mathcal{I}\mathcal{C}_x, j_{y\ast} \mathcal{C}_y[n])
= \op{dim}_{\Bbb{C}} \op{gMod}_{R-R} 
(B_x, \nabla_y
[n])
\end{displaymath} 
zeigen.
Nach \cite{SoBi}, Theorem 5.15 ist jedoch
$
\op{Mod}_{R-R} (B_x, \nabla_y)
$ 
graduiert frei als $R$-Rechtsmodul, 
und notieren wir $h^i_{y,x}$ die Zahl der im
Grad $i$ ben"otigten Erzeuger, 
so liefert Theorem 5.3 von loc.cit.\ in der Hecke-Algebra
die Formel
\begin{displaymath}
C^{\prime}_{x} = \sum_{y \in W} h^i_{y,x} q^{-1/2} \tilde{T}_y 
= \sum_{y \in W} h^i_{y,x}
v^i H_y
\end{displaymath}
in den Notationen von Lusztig bzw.\ von \cite{So-K}.
In anderen Worten haben wir $h^i_{y,x} = n^i_{y,x},$ und wegen 
$H^{\ast}_{B} (ByB/B) \cong R$
ergibt sich die behauptete Gleichheit der Dimensionen in jedem Grad.
Der zweite Fall folgt dual und damit ist das Lemma im Fall
$P=B$ vollst"andig bewiesen.
Im allgemeinen Fall folgt die Treuheit unseres Funktors ganz genauso,
aber f"ur den Dimensionsvergleich m"ussen wir uns noch etwas mehr anstrengen.
Wir behandeln hier nur die beiden F"alle $\mathcal{I} \mathcal{C} \rightarrow
\mathcal{I} \mathcal{C} $ und $\mathcal{I}\mathcal{C} \rightarrow \mathcal{C}$,
der verbleibende Fall kann dual behandelt werden.
Bezeichnet $\pi : G/B \twoheadrightarrow G/P$ die  Projektion, so haben
wir $\Bbb{H}_B \pi_\ast \mathcal{G} \cong \op{res}_{R-R}^{R-R^\iota} 
\Bbb{H}_B \mathcal{G}$
und $\mathcal{H}_B\pi^\ast \mathcal{F} \cong \Bbb{H}_B\mathcal{F}
\otimes_{R^\iota} R.$
Damit erhalten wir ein kommutatives Diagramm
\begin{displaymath}
\begin{array}{ccc}
\op{Der}_B (\mathcal{F}, \pi_\ast \mathcal{G} \left[ n\right]) & 
\overset{\sim}{\longrightarrow} & \op{Der}_B (\pi^\ast \mathcal{F}, 
\mathcal{G} \left[ n 
\right])\\
\downarrow & & \downarrow \\
\op{gMod}_{R-R^\iota} (\Bbb{H}_B \mathcal{F}, \Bbb{H}_B \pi_\ast 
\mathcal{G} \left[ n
\right]) & & \op{gMod}_{R-R}
(\Bbb{H}_B \pi^\ast \mathcal{F}, \Bbb{H}_B \mathcal{G} \left[ n
\right])\\
\parallel & &\parallel\\
\op{gMod}_{R-R^\iota} (\Bbb{H}_B \mathcal{F}, \op{res}_{R-R}^{R-R^\iota} 
\Bbb{H}_B \mathcal{G}
\left[ n \right]) &
\overset{\sim}{\longrightarrow} & \op{gMod}_{R-R} (\Bbb{H}_B \mathcal{F} 
\otimes_{R^\iota} R,
\Bbb{H}_B \mathcal{G} \left[ n \right])
\end{array}
\end{displaymath}
und mit der rechten oberen Vertikalen mu"s auch die linke 
Vertikale ein Isomorphismus
sein. Damit folgen die beiden F"alle $\mathcal{I}\mathcal{C} 
\rightarrow \mathcal{I}\mathcal{C}$
und $\mathcal{I}\mathcal{C} \rightarrow \mathcal{C}$ 
f"ur allgemeines $P$ aus dem Fall $P=B$.
\end{proof}

\section{Singul"are Bimoduln}
\begin{Bemerkung}
Sei $\mathcal{W}$ eine endliche Gruppe von 
Automorphismen eines endlichdimensionalen affinen Raums
$E$ "uber $\Bbb{Q}$, die von Spiegelungen 
erzeugt wird, und sei $\mathcal{S} \subset \mathcal{W}$
eine Wahl von einfachen Spiegelungen. 
Bezeichne $R$ die regul"aren Funktionen auf dem Raum
der Richtungsvektoren, graduiert durch die Vorschrift, 
da"s lineare Funktionen homogen
vom Grad 2 sein m"ogen.
So gibt es nach \cite{SoBi} bis auf Isomorphismus
 eindeutig bestimmte $\Bbb{Z}$-graduierte
$R$-Bimoduln $B_x = B_x (\mathcal{W}) = B_x (\mathcal{W},\mathcal{S}, E) 
\in R\op{-gMod-}R$ derart, da"s gilt
\begin{enumerate}
\item Die $B_x$ sind unzerlegbar.
\item
F"ur $e$ das neutrale Element ist $B_e=R.$
\item Ist $s \in \mathcal{S} $ eine einfache Spiegelung mit
$x s> x$, so gibt es eine Zerlegung
\begin{displaymath}
B_x \otimes_{R^s} R [1] \cong B_{xs} \oplus \bigoplus_{l(y) \leq l(x)}
{m(y)}B_y
\end{displaymath}
f"ur geeignete Vielfachheiten $m(y) \in \Bbb{N}$.
\end{enumerate}
Im "ubrigen besteht nach loc.cit.\
der Endomorphismenring dieser Bimoduln nur aus Skalaren,
insbesondere bleiben sie unzerlegbar unter Erweiterung der Skalare.
Sei nun $\mathcal{S}_{\iota} \subset \mathcal{S}$ eine 
Teilmenge der Menge der einfachen
Spiegelungen, $\mathcal{W}_\iota = \langle \mathcal{S}_\iota 
\rangle \subset \mathcal{W}$
ihr Erzeugnis, $w_\iota\in{\mathcal{W}}_\iota$ das l"angste 
Element  und $R^\iota$  der 
Teilring der ${\mathcal{W}_\iota}$-Invarianten. So behaupten wir
unter denselben Annahmen:
\end{Bemerkung}
\begin{Lemma}
F"ur jede Nebenklasse $\bar{x} \in \mathcal{W}/\mathcal{W}_\iota$ gibt es
genau einen unzerlegbaren $\Bbb{Z}$-graduierten $R$-$R^\iota$-Bimodul
$B_{\bar{x}}^\iota = 
B^\iota_{\bar{x}} (\mathcal{W},\mathcal{W}_\iota) \in R
\op{-gMod-}R^\iota$ mit der Eigenschaft, da"s f"ur $x \in \mathcal{W}$ den
l"angsten Repr"asentanten der Nebenklasse $\bar{x}$ gilt
\begin{displaymath}
\op{res}^{R-R^\iota}_{R-R} B_x \cong \bigoplus_{z \in 
\mathcal{W}_\iota} B^\iota_{\bar{x}}
[l(w_\iota)- 2l (z)]
\end{displaymath}
\end{Lemma}
\begin{proof}[Beweis]
Ohne Beschr"ankung der Allgemeinheit d"urfen wir annehmen, da"s 
$\mathcal{W}$ nur einen einzigen Fixpunkt hat.
Da es sich nach Annahme um eine rationale und mithin 
kristallographische Spiegelungsgruppe
handelt, finden wir dann $G \supset B \supset T$ eine komplexe halbeinfache
algebraische Gruppe $G$ mit Borel $B$ und maximalem 
Torus $T$ und Coxetersystem $(\mathcal{W},
\mathcal{S})$.
Identifizieren wir in $\mathcal{W}$-"aquivarianter 
Weise $H^2_{T} (\op{pt};\Bbb{Q})$ mit der 
homogenen Komponente $R^2$ von $R$, so gibt es 
nach \cite{So-L} einen Isomorphismus von
$\Bbb{Z}$-graduierten $R$-Bimoduln
\begin{displaymath}
B_x \cong \Bbb{H}_B \mathcal{I}\mathcal{C} (\overline{BxB / B})
\end{displaymath}
Ist $P=P_\iota$ 
eine Parabolische mit $G \supset P \supset B$, so zeigt der Zerlegungssatz
von \cite{BeLu}, angewandt auf die Projektion $p: G/B \twoheadrightarrow G/P,$
f"ur $x$ maximal in seiner $\mathcal{W}_{\iota}$-Nebenklasse schnell 
die Existenz einer Zerlegung
\begin{displaymath}
p_\ast \mathcal{I}\mathcal{C} (\overline{BxB/B}) \cong 
\bigoplus_{z\in \mathcal{W}_\iota} \mathcal{I}\mathcal{C}
(\overline{Bx P/P}) [l(w_\iota)-2l(z)]
\end{displaymath}
in $\op{Der}_B^+ (G/P)$. Andererseits haben wir aber 
$\Bbb{H}_B \circ p_\ast=  \op{res}^{R-R^\iota}_{R-R}   \circ\Bbb{H}_B$
und die $\Bbb{H}_B\mathcal{I}\mathcal{C}
(\overline{Bx P/P})$ sind unzerlegbar als graduierte $R$-$R^\iota$-Bimoduln,
da nach \ref{VTM} die Skalare ihre einzigen Endomorphismen vom
Grad $\leq 0$ sind.
\end{proof}
\section{Die Bimoduln zu Kippmoduln}
\begin{Bemerkung}
 Gegeben $y\in W$ bezeichne $ \hat{S}_{{y}}$ den Bimodul, der
  von links schlicht $\hat{S}$ ist, von rechts jedoch die mit ${y}$
  getwistete Operation $r 1_y = 1_y r^y$ von $\hat{S}$ tr"agt.
Gegeben ein Bimodul $B$ zu zwei kommutativen Ringen bezeichne
$\tilde{B}$ den Bimodul, der daraus durch Vertauschen der 
Linksoperation mit der Rechtsoperation entsteht.
\end{Bemerkung}
\begin{Satz}\label{BiKi}
Sei $\lambda \in \frak{h}^*_{\op{dom}}$ und seien $W_\lambda
\subset W_{\bar{\lambda}} \subset W$ wie in \ref{WLL}.
Gegeben $\bar{x} \in W_{\bar{\lambda}} / W_\lambda$
gilt f"ur die Deformation des unzerlegbaren Kippmoduls mit
h"ochstem Gewicht $w_{\bar{\lambda}} \bar{x} \cdot \lambda$
die Formel
\begin{displaymath}
\Bbb{V} K_{\hat{S}} (w_{\bar{\lambda}} \bar{x} \cdot \lambda)
\cong \widetilde{B}^\lambda_{\bar{x}} 
\otimes_S \hat{S}_{w_{\bar{\lambda}}}
\end{displaymath}
\end{Satz}
\begin{Bemerkung}
Unsere Bimoduln $B^\iota_{\bar{x}}$ sind graduiert frei von endlichem
Rang "uber $R$, folglich ist $\hat{R} \otimes_{R} B^\iota_{\bar{x}}$
schlicht der l"angs der Graduierung komplettierte Bimodul und er tr"agt
insbesondere eine Rechtsoperation von $\hat{R}^\iota$.
\end{Bemerkung}
\begin{proof}[Beweis]
Bekanntlich bleibt ein unzerlegbarer Kippmodul 
unzerlegbar, wenn wir ihn aus W"anden
r"ucken. Genauer gilt f"ur $\lambda, \mu \in \frak{h}^*_{\op{dom}}$ mit
$\lambda + X = \mu +X$ und $W_\mu =1$ und $x \in W_{\bar{\lambda}}$
maximal in seiner Nebenklasse $xW_{\lambda}$ notwendig
\begin{displaymath}
T^\mu_\lambda K (w_{\bar{\lambda}} x \cdot \lambda) =
K (w_{\bar{\mu}} x \cdot \mu)
\end{displaymath}
Da $T^\lambda_\mu T^\mu_\lambda$ eine Summe 
von $|W_\lambda|$ Kopien des Indentit"atsfunktors
ist, k"onnen wir uns beim Beweis des Satzes 
also auf den Fall $\lambda$ regul"ar 
beschr"anken.
Ist $x =st \ldots r$ eine reduzierte 
Darstellung durch einfache Spiegelungen von
$W_{\bar{\lambda}}$, so k"onnen wir 
dann $K_{\hat{S}} (w_{\bar{\lambda}}
x \cdot  \lambda)$
induktiv charakterisieren als den unzerlegbaren Summanden von
\begin{displaymath}
\vartheta_r \ldots \vartheta_t \vartheta_s 
\Delta_{\hat{S}} (w_{\bar{\lambda}}
\cdot \lambda)
\end{displaymath}
der nicht isomorph ist zu einem 
$K_{\hat{S}} (w_{\bar{\lambda}}y \cdot \lambda)$
f"ur $y < x$.
Wenden wir $\Bbb{V}$ an, so erhalten wir 
daraus den unzerlegbaren Summanden von
\begin{displaymath}
\hat{S} \otimes_{\hat{S}^r} 
\hat{S} \ldots \otimes_{\hat{S}^t} \hat{S}
\otimes_{\hat{S}^s} \hat{S}_{w_{\bar{\lambda}}}
\end{displaymath}
der nicht ``schon vorher vorkam''.
Das ist aber nach der Definition der speziellen Bimoduln  genau 
$\widetilde{B}_x \otimes_S \hat{S}_{w_{\bar{\lambda}}}$.
\end{proof}
\section{Erg"anzungen zum Wechsel der Gruppe}
\begin{Bemerkung}
F"ur $G$ eine komplexe zusammenh"angende algebraische Gruppe 
setzen wir $A_G=H^\ast(BG).$ 
Ist $X$ eine komplexe 
algebraische $G$-Variet"at und sind $\mathcal{F}, \mathcal{G} \in
\op{Der}^+_{G} (X)$ Objekte der "aquivarianten derivierten Kategorie,
so bilden wir den graduierten $A_G$-Modul
$$\op{Der}_G (\mathcal{F}, \mathcal{G} [\ast])
=\bigoplus_n \op{Der}_G (\mathcal{F}, \mathcal{G} [n])$$
  \end{Bemerkung}
\begin{Proposition}\label{GruWe}
Seien $G \supset H$ eine zusammenh"angende komplexe algebraische Gruppe und
eine zusammenh"angende abgeschlossenen Untergruppe.
Sei $X$ eine algebraische $G$-Variet"at und seien 
$\mathcal{F}, \mathcal{G} \in
\op{Der}_{G} (X)$ konstruierbare Komplexe.
Ist $\op{Der}_G (\mathcal{F}, \mathcal{G} [\ast])$ graduiert frei
"uber $A_G,$
so induziert die offensichtliche Abbildung eine Bijektion
\begin{displaymath}
A_H \otimes_{A_G} 
 \op{Der}_G (\mathcal{F},\mathcal{G} 
[\ast]) \overset{\sim}{\rightarrow} 
\op{Der}_H (\mathcal{F},\mathcal{G} [\ast])
\end{displaymath}
\end{Proposition}
\begin{proof}[Beweis]
Wir betrachten die konstante Abbildung
$k: X \ra \op{pt}$ und den volltreuen Funktor
$
\gamma_G : \op{Der}_G^{\op{c}} (\op{pt}) \ra A_G\op{-dgDer} 
$
nach \cite{BeLu} und beachten
$$
\op{Der}_G (\mathcal{F}, \mathcal{G} [\ast]) = H^\ast \gamma_G k_{\ast}
\op{Hom}(\mathcal{F},\mathcal{G})
$$
wo wir $\op{Hom}(\mathcal{F}, \mathcal{G})$ 
in $\op{Der}_G^+ (X)$ bilden und
$k_{\ast}$ das direkte Bild in $\op{Der}_G^+ (\op{pt})$ meint.
Ist das nun ein freier $A_G$-Modul,
so ist $\gamma_G k_\ast \op{Hom} (\mathcal{F},\mathcal{G})$ bereits
quasiisomorph zu seiner Kohomologie und diese 
Kohomologie ist homotopieprojektiv in 
$A_G\op{-dgMod}.$
Mit dem derivierten Funktor
$A_H \otimes^{L}_{A_G}: A_G\op{-dgDer}\ra A_H\op{-dgDer}$
 haben wir nach \cite{BeLu}, 12.7.1 weiter kanonisch
\begin{displaymath}
A_H \otimes^{L}_{A_G} \circ\gamma_G 
= \gamma_{H} \circ \op{res}^H_G
\end{displaymath}
und bei homotopieprojektiven Objekten
$M \in A_G\op{-dgMod}$ haben wir zus"atzlich
\begin{displaymath}
A_H \otimes^L_{A_G} M
= A_H \otimes_{A_G} M
\end{displaymath}
Da $k_{\ast}$ und $\op{Hom}$ mit der 
Restriktion der Gruppenoperation vertauschen,
zeigt das die Proposition.
\end{proof}

\section{Geometrie  der Filtrierung}
\begin{Bemerkung}
Gegeben ein Ring $R,$  ein Ringhomomorphismus
$R \ra \Bbb{C} [[t]],$ 
drei $R$-Moduln $H,H',H''$ und eine $R$-bilineare Abbildung
\begin{displaymath}
\varphi : H \times H^\prime \rightarrow H^{''}
\end{displaymath}
k"onnen wir auf $\Bbb{C} [[t]]\otimes_R H$ eine Filtrierung
erkl"aren durch die Vorschrift
\begin{displaymath}
F^i (\Bbb{C} [[t]] \otimes_R H) = \left\{
h \left| \begin{array}{l}
\varphi (h,h^\prime)\in t^i \Bbb{C} [[t]] \otimes_R H^{''}\\
\text{f"ur alle } h^\prime \in \Bbb{C} [[t]] \otimes_R H^\prime
\end{array} \right. \right\}
\end{displaymath}
und erhalten dann nat"urlich auch 
eine induzierte Filtrierung auf $\Bbb{C} \otimes_R H,$
deren Subquotienten wir $\bar{F}^i(\Bbb{C} \otimes_R H)$ notieren.
\end{Bemerkung}
\begin{Bemerkung}\label{NZ}
Ist zum Beispiel $\hat{S}$ die Komplettierung von 
$S = \mathcal{O} (\frak{h}^*)$
l"angs der nat"urlichen Graduierung und $\hat{S} 
\twoheadrightarrow \Bbb{C} [[t]]$ die Restriktion
auf die Gerade $t\rho$ wie in \ref{DAF}, 
so liefert die durch Komposition gegebene Paarung
\begin{displaymath}
\begin{array}{r}
\op{Hom}_{\frak{g} - \hat{S}} 
(\Delta_{\hat{S}} (\lambda), K_{\hat{S}} (\mu)) \times 
\op{Hom}_{\frak{g}-\hat{S}} (K_{\hat{S}} (\mu), 
\nabla_{\hat{S}} (\lambda))\\[2mm] \rightarrow
\op{Hom}_{\frak{g}-\hat{S}} (\Delta_{\hat{S}} (\lambda), 
\nabla_{\hat{S}} (\lambda))
\end{array}
\end{displaymath}
die Andersen-Filtrierung auf dem Raum
$\op{Hom}_{\frak{g}} (\Delta (\lambda), K (\mu)),$
der ja durch \ref{BCH} mit 
$
\op{Hom}_{\frak{g}-\hat{S}} (\Delta_{\hat{S}} (\lambda), 
K_{\hat{S}} (\mu)) \otimes_{\hat{S}}
\Bbb{C}$ identifiziert werden kann.
Nat"urlich gibt es auch ein 
$p\in \hat{S}$ derart, da"s $\Bbb{V}$ einen Isomorphismus
\begin{displaymath}
\op{Hom}_{\frak{g}-\hat{S}} (\Delta_{\hat{S}} (\lambda), 
\nabla_{\hat{S}} (\lambda))
\overset{\sim}{\rightarrow} p \op{Hom}_{\hat{S}} 
(\Bbb{V} \Delta_{\hat{S} }(\lambda),
\Bbb{V} \nabla_{\hat{S}} (\lambda))
\end{displaymath}
induziert, und mit diesem $p$ k"onnen wir unsere 
Paarung umschreiben zu der wieder durch
Komposition gegebenen Paarung
\begin{displaymath}
\begin{array}{r}
\op{Hom}_{\hat{S}^{\lambda} -\hat{S}} 
(\Bbb{V} \Delta_{\hat{S}} (\lambda),
\Bbb{V} K_{\hat{S}}(\mu)) \times  
\op{Hom}_{\hat{S}^{\lambda} - \hat{S}}
(\Bbb{V} K_{\hat{S}} (\mu), \Bbb{V} \nabla_{\hat{S}} (\lambda))\\[2mm]
\rightarrow p \op{Hom}_{\hat{S}^{\lambda}-\hat{S}} 
(\Bbb{V} \Delta_{\hat{S}}
(\lambda), \Bbb{V} \nabla_{\hat{S}} (\lambda))
\end{array}
\end{displaymath}
Hier kann ein m"ogliches $p$ dadurch bestimmt werden, da"s 
unsere Paarung f"ur 
$\lambda = \mu$ eine surjektive Abbildung
liefern mu"s.
Nun wechseln wir die Parameter, w"ahlen $\lambda \in \frak{h}^*_{\op{dom}}$ und
notieren $\lambda_{\bar{x}} = w_{\bar{\lambda}} \bar{x} \cdot \lambda$
f"ur $\bar{x} \in W_{\bar{\lambda}}/ W_{\lambda}$.
Um uns nicht hoffnungslos in der Notation zu verheddern, 
definieren wir weiter eine
Variante $\widetilde{\Bbb{V}}$ von $\Bbb{V}$ durch die Vorschrift
\begin{displaymath}
\widetilde{\Bbb{V}} M = \hat{S}_{w_{\bar{\lambda}}} 
\otimes_{\hat{S}}
\widetilde{\Bbb{V}M}
\end{displaymath}
so da"s sich \ref{BiKi} vereinfacht zu
$
\widetilde{\Bbb{V}} K_{\hat{S}} (\lambda_{\bar{x}}) 
\cong \hat{B}^\lambda_{\bar{x}}
,$
wobei der Hut die Komplettierung l"angs der Graduierung meint.
Mit weniger M"uhe pr"uft man auch $$\widetilde{\Bbb{V}} 
\Delta_{\hat{S}} (\lambda_{\bar{y}})
\cong \widetilde{\Bbb{V}} \nabla_{\hat{S}} 
(\lambda_{\bar{y}}) \cong \hat{S}^\lambda_{\bar{y}}$$
in $\hat{S}\op{-mod-}\hat{S}^\lambda$, 
womit wieder der Bimodul gemeint ist, der von links
schlicht $\hat{S}$ ist, von rechts jedoch 
die mit $\bar{y}$ getwistete Operation
$r 1_y = 1_y r^y$ von $\hat{S}^\lambda$ tr"agt.
Ersetzen wir also $\Bbb{V}$ durch $\widetilde{\Bbb{V}}$, 
so landen wir bis auf einen
Twist der $\hat{S}$-Struktur um $w_{\bar{\lambda}}$ bei der Paarung
\begin{displaymath}
  \begin{array}{r}
\op{Hom}_{\hat{S}-\hat{S}^\lambda} (\hat{S}^\lambda_{\bar{y}},
\hat{B}^\lambda_{\bar{x}})\times \op{Hom}_{\hat{S}-\hat{S}^\lambda}
(\hat{B}^\lambda_{\bar{x}}, \hat{S}^\lambda_{\bar{y}})
\\[2mm] 
\rightarrow p_1 
\op{Hom}_{\hat{S}-\hat{S}^\lambda}({\hat{S}^\lambda_{\bar{y}}}, 
\hat{S}^\lambda_{\bar{y}})
\end{array}
\end{displaymath}
von $\hat{S}$-Moduln und unsere Filtrierung 
entspricht eben der Filtrierung, die sich
hier ergibt, wenn wir unser
$\hat{S} \twoheadrightarrow \Bbb{C} [[t]]$ dadurch
ab"andern, da"s wir es mit $w_{\bar{\lambda}}$
vertwisten. Hier meint $p_1$ das Bild
von $p$ unter  $w_{\bar{\lambda}}.$
Da es bei der Wahl von $p_1$ eh nicht auf
Einheiten von $\hat{S}$ ankommt,  k"onnen wir auch $p_1$ bereits im
noch nicht komplettierten Ring w"ahlen, und
die entsprechende Paarung ``vor Komplettierung'' 
\begin{displaymath}
  \begin{array}{r}
\op{Hom}_{S-S^\lambda} (S^\lambda_{\bar{y}}, 
B^\lambda_{\bar{x}}) \times \op{Hom}_{S-S^\lambda}
(B^\lambda_{\bar{x}}, S^\lambda_{\bar{y}}) \\[2mm] \rightarrow
p_1 \op{Hom}_{S-S^\lambda} (S^\lambda_{\bar{y}}, S^\lambda_{\bar{y}})
\end{array}
\end{displaymath}
liefert am Ende auch denselben filtrierten $\Bbb{C}$-Vektorraum. 
Diese Paarung hinwiederum interpretieren wir
nun geometrisch.
\end{Bemerkung}
  
\begin{Bemerkung}
Bezeichnet $X \subset \frak{h}^*$ das Gitter der 
ganzen Gewichte, so gibt es 
ein Paar $G^\vee \supset T^\vee$ bestehend aus 
einer reduktiven zusammenh"angenden komplexen
algebraischen Gruppe mit einem maximalen Torus derart, da"s $X = X(T^\vee)$ 
die Gruppe seiner Einparameteruntergruppen ist 
und da"s f"ur ihre  Weylgruppe gilt $W( G^\vee, T^\vee)=W_{\bar{\lambda}}.$ 
In $G^\vee$ w"ahlen wir dann eine Borel $B^\vee$ zu $\mathcal{S} 
\cap W_{\bar{\lambda}}$ und eine Parabolische
$P^\vee \supset B^\vee$ zu $W_\lambda$.
Bezeichnet nun $$\mathcal{I} \mathcal{C}_{\bar{x}} 
= \mathcal{I}\mathcal{C}
(\overline{B^\vee \bar{x} P^\vee/P^\vee})$$
den Schnittkohomologiekomplex der entsprechenden Schubertvariet"at
und $\mathcal{C}_{\bar{y}}$ die konstante 
perverse Garbe auf $B^\vee\bar{y} P^\vee/P^\vee$,
so haben wir $B^\lambda_{\bar{x}} \cong 
\Bbb{H}_{B^\vee} \mathcal{I}\mathcal{C}_{\bar{x}}$
und unsere Paarung ``vor Komplettierung'' vom Schlu"s der vorhergehenden
Bemerkung
l"a"st sich mithilfe von \ref{VTM} interpretieren als die
immer noch durch Komposition gegebene Paarung
\begin{displaymath}
  \begin{array}{r}
\op{Der}^\ast_{B^\vee} (j_{\bar{y}!}\mathcal{C}_{\bar{y}}, 
\mathcal{I}\mathcal{C}_{\bar{x}})
\times \op{Der}^\ast_{B^\vee} (\mathcal{I}\mathcal{C}_{\bar{x}}, 
j_{\bar{y}\ast} \mathcal{C}_{\bar{y}})\\[2mm]
\rightarrow \op{Der}^\ast_{B^\vee} (j_{\bar{y}!} 
\mathcal{C}_{\bar{y}}, j_{\bar{y}\ast}
\mathcal{C}_{\bar{y}})
\end{array}
\end{displaymath}
In der Tat mu"s hier auf der rechten Seite kein $p$-Faktor 
erg"anzt werden, da der R"uckzug auf die dicke Zelle rasch zeigt,
da"s  im Fall $\bar{x}=\bar{y}$ unsere Paarung eine Surjektion liefert.
Die Frage ist also, welchen filtrierten Vektorraum diese Paarung von
$A_{B^\vee}$-Moduln liefert unter dem
Homomorphismus $A_{B^\vee} \twoheadrightarrow \Bbb{C} [t],$
der durch die Einbettung $\Bbb{C}^\times \hookrightarrow T^\vee$ 
mit Parameter $w_{\bar{\lambda}}
\rho$ gegeben wird.
Nach \ref{GruWe} f"uhrt uns jedoch die Spezialisierung auf 
$\Bbb{C} [t]$ zur durch Komposition gegebenen Paarung
\begin{displaymath}
  \begin{array}{r}
\op{Der}^\ast_{\Bbb{C}^\times} (j_{\bar{y}!}\mathcal{C}_{\bar{y}}, 
\mathcal{I}
\mathcal{C}_{\bar{x}}) \times \op{Der}^\ast_{\Bbb{C}^\times}
(\mathcal{I}\mathcal{C}_{\bar{x}}, j_{\bar{y}\ast} 
\mathcal{C}_{\bar{y}})\\[2mm]
\rightarrow \op{Der}^\ast_{\Bbb{C}^\times} (j_{\bar{y}!}
\mathcal{C}_{\bar{y}},
j_{\bar{y}\ast} \mathcal{C}_{\bar{y}})
\end{array}
\end{displaymath}
Bezeichne $\bar{y}$ im Folgenden auch den Punkt 
$\bar{y}P^\vee$ von $G^\vee/P^\vee.$
F"ur ein geeignetes Produkt $U$ von Wurzelgruppen aus $G^\vee$ definiert die
Multiplikation $u \mapsto u{\bar{y}}$ eine Einbettung 
$U \hookrightarrow G^\vee / P^\vee,$ 
deren Bild eine zur Zelle $B^\vee{\bar{y}} P^\vee / P^\vee$ 
transversale Zelle  ist, die unter unserem
$\Bbb{C}^\times$ auf  $\bar{y}$ kontrahiert wird.
Setzen wir $Z = U{\bar{y} }\cap \overline{B^\vee\bar{x} P^\vee/P^\vee},$
so wird auch $Z$ von $\Bbb{C}^\times$ auf $\bar{y}$ kontrahiert, und
bezeichnet $a : Z \hookrightarrow G^\vee/P^\vee$ die Einbettung,
so wird das Restringieren auf $Z$ unsere Paarung von
eben nicht "andern.
Setzen wir nun  $d =\dim B^\vee{\bar{y} }P^\vee/P^\vee$ und bezeichnen mit
$j:\op{pt}\hra Z$ die Einbettung von $\bar{y}$ 
und mit $\underline{\op{pt}}$ die konstante Garbe auf einem Punkt, so 
erhalten wir 
$
a^* j_{\bar{y}!} \mathcal{C}_{\bar{y}} \cong j_* \underline{\op{pt}}
[d] \cong a^*j_{\bar{y}_*} \mathcal{C}_{\bar{y}}
$ 
und $a^* \mathcal{I}\mathcal{C}_{\bar{x}} \cong \mathcal{I} \mathcal{C}
 [d]$ wird der verschobene Schnittkohomologiekomplex $\mathcal{I} \mathcal{C}
=\mathcal{I} \mathcal{C}(Z)$ von $Z$ 
und unsere Paarung verwandelt  sich in die durch Komposition gegebene
Paarung
\begin{displaymath}
\begin{array}{r}
\op{Der}^*_{\Bbb{C}^\times} (j_* \underline{\op{pt}}, 
\mathcal{I}\mathcal{C}) \times \op{Der}^*_{\Bbb{C}^\times}
(\mathcal{I}\mathcal{C} , j_* \underline{\op{pt}})\\[2mm]
 \rightarrow  \op{Der}^*_{\Bbb{C}^\times} (j_* \underline{\op{pt}}, j_* \underline{\op{pt}})
\end{array}
\end{displaymath}
Nun k"onnen wir den ersten unserer gepaarten Moduln 
identifizieren mit dem Kohalm $j^! \mathcal{I}
\mathcal{C} $ und den Zweiten mit dem Dualen 
des Halms $j^* \mathcal{I}\mathcal{C} $ und
unsere Paarung liefert folglich denselben 
filtrierten Vektorraum wie die Einbettung
von freien $\Bbb{C}[t]$-Moduln 
$j^! \mathcal{I}\mathcal{C}  \hookrightarrow j^* \mathcal{I}
\mathcal{C} .$
Jetzt besagt aber das ``Fundamental example'' 
aus Abschnitt 14 von \cite{BeLu}  gerade, da"s der Kokern dieser Einbettung
identifiziert werden kann mit der 
Schnittkohomologie der projektiven Variet"at 
$\overline{Z} = (Z- \{\overline{y}\})
/\Bbb{C}^\times$ verschoben um Eins, 
in Formeln also mit $IC (\overline{Z}) [1]$, 
aufgefa"st als $\Bbb{C} [t]$-Modul
aufzufassen ist in der Weise, 
da"s $t$ als der Lefschetzoperator wirkt. Der harte Lefschetz 
f"ur die Schnittkohomologie aus \cite{BBD} sagt uns nun,
da"s die fragliche Filtrierung auf 
$\Bbb{C} \otimes_{\Bbb{C}[t]} j^! \mathcal{I}\mathcal{C} $
"ubereinstimmt mit der durch die 
$\Bbb{Z}$-Graduierung gegebenen Filtrierung, und diese kann
bekanntlich durch Kazhdan-Lusztig-Polynome beschrieben werden.
Genauer erhalten wir
\begin{displaymath}
\begin{array}{ccl}
\bar{F}^i (\Bbb{C} \otimes_{\Bbb{C}[t]} j^! \mathcal{I}\mathcal{C}) &
\cong & H^{-i} (j^! \mathcal{I}\mathcal{C})\\
& \cong & \op{Der} (\mathcal{C}_{\bar{y}} [i],j^!_{\bar{y}} 
\mathcal{I}
\mathcal{C}_{\bar{x}})\\
&\cong & \op{Der} (j_{\bar{y}!} \mathcal{C}_{\bar{y}} [i], 
\mathcal{I}
\mathcal{C}_{\bar{x}})
\end{array}
\end{displaymath}
und dieser Raum hat bekanntlich die Dimension $h^i_{y,x}$ f"ur $y,x$ 
die l"angsten
Repr"asentanten von $\bar{y}, \bar{x}$.
Das ist also auch die Dimension des $i$-ten Subquotienten der 
Andersen-Filtrierung auf
\begin{displaymath}
\op{Hom}_{\frak{g}}(\Delta (\lambda_{\bar{y}}), 
K(\lambda_{\bar{x}}))
\end{displaymath}
\end{Bemerkung}

\providecommand{\bysame}{\leavevmode\hbox to3em{\hrulefill}\thinspace}
\providecommand{\MR}{\relax\ifhmode\unskip\space\fi MR }
\providecommand{\MRhref}[2]{%
  \href{http://www.ams.org/mathscinet-getitem?mr=#1}{#2}
}
\providecommand{\href}[2]{#2}

\end{document}